\documentclass[11pt,a4paper]{article}

\usepackage{amsmath,amssymb, bbm}
\usepackage{color}
\newcommand{\R}{\mathbb{R}}

\newcommand{\Z}{\mathbb{Z}}

\def\cA{{\mathcal A}}
\def\cC{{\mathcal C}}

\def\cS{{\mathcal S}}

\newcommand{\ee}{\varepsilon}
\renewcommand{\aa}{\alpha}
\renewcommand{\div}{{\rm div}\,}

\newcommand{\Sum}{\displaystyle \sum}

\def\d{\partial}
\def\ddl{\dot \Delta_l}
\def\ddj{\dot \Delta_j}
\def\ddq{\dot \Delta_q}

\def\tilde{\widetilde}
\def\hat{\widehat}
\newcommand{\D}{\Delta}

\newcommand{\n}{\nabla}

\newcommand{\du}{\delta u}

\newcommand{\mub}{\bar{\mu}}
\newcommand{\kb}{\bar{\kappa}}
\newcommand{\trho}{\tilde{\rho}}
\newcommand{\tu}{\tilde{u}}
\newcommand{\tP}{\tilde{P}}
\newcommand{\tv}{\tilde{v}}
\newcommand{\tQ}{\tilde{Q}}
\newcommand{\kmu}{\frac{\kb}{\mub^2}}
\newcommand{\ub}{\bar{u}}
\newcommand{\vb}{\bar{v}}
\newcommand{\wb}{\bar{w}}
\newcommand{\rhob}{\bar{\rho}}
\newcommand{\ro}{\rho_0}
\newcommand{\Pb}{\bar{P}}
\newcommand{\Qb}{\bar{Q}}
\newcommand{\Au}{A_{\ub}}
\newcommand{\Av}{A_{\vb}}
\newcommand{\adjv}{\mbox{adj}(D X_{\vb})}
\newcommand{\adju}{\mbox{adj}(D X_{\ub})}

\newcommand{\Jv}{J_{\vb}}
\newcommand{\dub}{\delta \ub}
\newcommand{\dPb}{\delta \Pb}
\newcommand{\dwb}{\delta \wb}
\newcommand{\dQb}{\delta \Qb}
\newcommand{\rhok}{\rho_{\kb}}
\newcommand{\uk}{\bar{u}_{\kb}}
\newcommand{\Pk}{\bar{P}_{\kb}}
\newcommand{\dP}{\delta P}

\newtheorem{thm}{Theorem}

\newtheorem{prop}{Proposition}

\newtheorem{rem}{Remark}

\title{Lagragian methods for a general inhomogeneous incompressible Navier-Stokes-Korteweg system with variable capillarity and viscosity coefficients}

\author{Cosmin Burtea\footnote{Universit\'e Paris-Est Cr\'eteil, Laboratoire d'Analyse et de Math\'ematiques Appliqu\'ees (UMR 8050), 61 Avenue du G\'en\'eral de Gaulle, 94 010 Cr\'eteil Cedex (France). E-mail: cosmin.burtea@u-pec.fr}, Fr\'ed\'eric Charve\footnote{Universit\'e Paris-Est Cr\'eteil, E-mail: frederic.charve@u-pec.fr}}

\date{}
\begin{document}

\maketitle

\begin{abstract}
We study the inhomogeneous incompressible Navier-Stokes system endowed with a general capillary term. Thanks to recent methods based on Lagrangian change of variables, we obtain local well-posedness in critical Besov spaces (even if the integration index $p\neq 2$) and for variable viscosity and capillary terms. In the case of constant coefficients and for initial data that are perturbations of a constant state, we are able to prove that the lifespan goes to infinity as the capillary coefficient goes to zero, connecting our result to the global existence result obtained by Danchin and Mucha for the incompressible Navier-Stokes system with constant coefficients.
\end{abstract}



\section{Introduction}

\subsection{Presentation of the model}

A classical model to study the dynamics of a fluid endowed with internal capillarity (in the diffuse interface setting) is the following general compressible Korteweg system:
\begin{equation}
\begin{cases}
&\d_t\rho+\div (\rho u)=0,\\
&\d_t (\rho u)+\div (\rho u\otimes u)-\cA u+\nabla(P(\rho))=\div K +\rho f,\\
\end{cases}
\end{equation}
where $\rho(t,x)\in \R$ and $u(t,x)\in \R^n$ respectively denote the density and velocity of the fluid at the position $(t,x)\in [0,\infty[ \times \R^n$ with $n\geq 2$.
The diffusion operator is $\cA u =\div(\mu(\rho) D(u)) +\n (\lambda(\rho) \div u)$, and the differential operator $D$ is defined by $D(u)=\n u+Du$ (i.-e. $D(u)_{i,j}= \d_i u^j+ \d_j u^i$, we chose here the double of the classical symmetric gradient).\\
The scalar functions $\mu$ and $k$ are smooth functions $\R \rightarrow \R$, the pressure $P$ is a given function (usually chosen as the Van der Waals pressure) and the general capillary tensor is written as:
$$
\div K = \nabla\left(\rho k(\rho) \Delta \rho +\frac{1}{2} \Big(k(\rho) + \rho k'(\rho)\Big) |\n \rho|^2\right) -\div\Big(k(\rho) \n \rho \otimes \n \rho\Big).
$$
In this article we are interested in the incompressible counterpart, namely the incompressible inhomogeneous capillary Navier-Stokes system:
\begin{equation}
\begin{cases}
\d_t \rho +u\cdot \n \rho=0,\\
\rho (\d_t u +u\cdot \n u) -\mub  \div \left(\mu (\rho) D(u)\right) +\n P= -\kb \div \left(k(\rho) \n \rho \otimes \n \rho\right),\\
\div u=0,\\
(\rho,u)_{|t=0}=(\rho_0, u_0),
\end{cases}
\label{NSIK}
\tag{$NSIK$}
\end{equation}
where the unknown are now $(\rho, u, P)$. We focus on solutions around the stable equilibrium state $(\bar{\rho},0)$ and in order to simplify the notations we will choose $\bar{\rho}=1$ and we introduce the parameters $\mub$ and $\bar{\kappa}$ in order to normalise the variable coefficients with $\mu(1)=k(1)=1$.
\begin{rem}
\sl{Due to the incompressibility condition the usual second viscous term, namely $\n (\lambda(\rho) \div u)$, vanishes in the pressure term and for the same reason we directly write the capillary term as a general divergence form (the case of a constant function $k\equiv 1$ corresponds to the classical capillary coefficients $\rho \n \Delta \rho$ and $-\n \rho \Delta \rho$ which are equivalent as their difference is a gradient, absorbed by the pressure term).}
\end{rem}

Let us first give a few words about the classical compressible Navier-Stokes-Korteweg system in the case of a \emph{constant coefficient} $k(\rho)$:
\begin{equation}
\begin{cases}
&\d_t\rho+\div (\rho u)=0,\\
&\d_t (\rho u)+\div (\rho u\otimes u)-\cA u+\nabla(P(\rho))=k\rho\nabla\Delta \rho.\\
\end{cases}
\tag{NSK}
\label{NSK}
\end{equation}
As pointed out by Danchin and Desjardins in \cite{DD} there is a strong coupling between the compressible part of the velocity and the density fluctuation $\ro-1$ that helps regularizing $\ro-1$ (the incompressible part of the velocity being totally decoupled with it). This is why $\ro-1$ features a parabolic regularization for all frequencies, unlike for the classical compressible Navier-Stokes system where this occurs only for low frequencies. More precisely they prove that if $\ro-1\in \dot{B}_{2,1}^\frac{n}{2}$ and $\ro$ bounded from below by a positive constant, and if $u_0 \in \dot{B}_{2,1}^{\frac{n}{2}-1}$ the system has a unique local solution $(\rho,u)$ on $[0,T]$ with
$$
\begin{cases}
\rho-1 \in C_T \dot{B}_{2,1}^\frac{n}{2} \cap L_T^1 \dot{B}_{2,1}^{\frac{n}{2}+2}\\
u \in C_T \dot{B}_{2,1}^{\frac{n}{2}-1} \cap L_T^1 \dot{B}_{2,1}^{\frac{n}{2}+1},
\end{cases}
$$
And in the case of small data (with more regular initial density fluctuation $\ro-1 \in \dot{B}_{2,1}^{\frac{n}{2}-1} \cap \dot{B}_{2,1}^\frac{n}{2}$), near a stable constant state ($P'(1)>0$) they obtain a global solution.

Let us also mention the works of Haspot (see \cite{Has1, Has2, Has4, Has5}) where the Korteweg system is studied in general settings, or with minimal assumptions on the initial data (or Besov indices) in the case of particular viscosity or capillarity coefficients.

Another way of modelling the capillarity (still in the diffuse interface setting) consists in considering a non-local capillary term, featuring only one derivative. This was first suggested by Van der Waals and re-discovered by Rohde (see \cite{Ro2, Rohdeorder}). As studied in \cite{Has2, CH1, CH4, FCorder} for small perturbations of a stable constant case (in critical spaces for integrability index $p=2$) the regularity structure of the density fluctuation is closer to the classical compressible Navier-Stokes case: parabolic regularization in low frequencies, damping in the high frequency regime (the threshold depending on $\mu^2/\kappa$). Moreover there is also a strong coupling between the density and the gradient part of the velocity, and in the cited works is also studied the transition from the non-local to the local capillary models. We also mention \cite{FCLarge} for the same study without smallness and stability assumptions.

In the case of the inhomogeneous incompressible Korteweg, the velocity is divergence-free so that no part of it can combine with the density fluctuation and improve its regularization. Any energy method with use of symmetrizers is completely useless. This system is much less studied than the compressible version. We can for example refer to \cite{B2, YYL} both of them in the case of constant viscosities and constant capillary coefficient $\n \rho \Delta \rho$. The first article provides local solutions for bounded domains and regular initial data in Sobolev spaces $W^{s,p}$ ($p\geq 1$) and requires $L^2$-assumptions. The second article requires regular initial data in energy spaces ($H^3$) and obtains local solutions as well as convergence towards the Euler system as $\kappa$ and $\mu$ go to zero. It is crucial in their article to take advantage of the $L^2$-structure to deal with the capillary term as $(\nabla \rho \cdot \Delta \rho|u)=(u\cdot \nabla \rho|\Delta \rho)= -(\d_t \rho|\Delta \rho)$.

In the present article, we wish to obtain well-posedness results for initial data in critical spaces (minimal regularity assumptions), for general variable capillary and viscosity coefficients and for general integrability index ($p$ not necessarily equal to $2$) and we will follow Lagrangian methods developped by \cite{RP2, Dlagrangien2} (and now frequently used, we refer for example to \cite{Noboru, PZ1, PZ2}). Let us precise that this method was first introduced by Hmidi in \cite{TH1} and extended by Danchin in \cite{Dlagrangien}. The incompressibility condition in our system leads us to adapt the arguments of \cite{RP2} but the additional capillary term will force us to mix them with arguments from \cite{Dlagrangien2} that extends the methods from \cite{RP2} in the compressible setting and therefore deals with additional external force terms that are only bounded in time (which imposes bounds for the lifespan).

\begin{rem}
\sl{We present here results for general capillary and viscosity terms in the incompressible setting but various non-local capillary terms can also be considered in compressible or quasi-incompressible cases (see for example \cite{Abels, Aki}).}
\end{rem}

\subsection{Statement of the results}
The main goal of this article is to state results for \emph{general smooth viscosity and capillary coefficients}, \emph{without any smallness assumptions} on the initial density fluctuation $\rho_0-1$, and in \emph{truly critical spaces}, that is without any extra low frequencies assumptions, all of this in a general $L^p$-setting (without energy methods). Thanks to the Lagrangian methods we will be able to give a very short proof of these results.

For this we will use the new class of estimates obtained by the first author in \cite{CB} instead of the classical maximal regularity estimates obtained by Danchin and Mucha in \cite{RP2, Dlagrangien2}. As in \cite{CB} our result will be given for $n=2$ with $p\in (1,4)$ and $n=3$ with $p\in(6/5,4)$.
\\

If we wish to recover the full usual range for $n,p$ ($p\in[1,2n[$) we need to assume as in \cite{RP2, Dlagrangien} that $\rho_0-1$ is small. Like in \cite{RP2} we obtain a local solution. Finally if we also assume that the initial velocity is small, although we are unable to obtain global existence as in \cite{RP2}, we are able to prove that the lifespan can be bounded from below by $\frac{C}{\kb}$, which goes to infinity as $\kb$ goes to zero. Moreover the solution converges on any time interval $[0,T]$ to the global solution of the inhomogeneous incompressible Navier-Stokes system given by \cite{RP2}.
\\

Let us now state in detail these results:

\begin{thm}
\sl{Assume that $n=2$ and $p\in(1,4)$ or $n=3$ and $p\in(\frac{6}{5},4)$. Let the initial data $(\rho_0, u_0)$ satisfy:
$$
\div u_0=0, \mbox{ }\inf_{x\in \R^n}\rho_0(x)>0, \mbox{ }u_0\in \dot{B}_{p,1}^{\frac{n}{p}-1}, \mbox{ } \rho_0-1, \mbox{ }\nabla \rho_0 \in \dot{B}_{p,1}^\frac{n}{p}.
$$
Then there exists a positive time $T>0$ and a unique solution $(\rho, u, \nabla P)$ of \eqref{NSIK} with
$$
\begin{cases}
\rho-1, \mbox{ }\nabla \rho \in C_T(\dot{B}_{p,1}^\frac{n}{p}), \quad \d_t \rho \in L_T^\infty \dot{B}_{p,1}^{\frac{n}{p}-1} \cap L_T^2 \dot{B}_{p,1}^\frac{n}{p},\\
u\in C_T(\dot{B}_{p,1}^{\frac{n}{p}-1}), \quad (\d_t u, \nabla^2 u, \nabla P) \in L_T^1(\dot{B}_{p,1}^{\frac{n}{p}-1}).
\end{cases}
$$
}
\label{T1}
\end{thm}
\begin{rem}
\sl{At the end of the article we give more details about the dependency of $T$ in terms of the parameters}.
\end{rem}

\begin{thm}
\sl{Assume that $n\geq 2$ and $p\in[1,2n)$. Let the initial data $(\rho_0, u_0)$ satisfy:
$$
\div u_0=0, \mbox{ }\inf_{x\in \R^n}\rho_0(x)>0, \mbox{ }u_0\in \dot{B}_{p,1}^{\frac{n}{p}-1}, \mbox{ } \rho_0-1, \mbox{ }\nabla \rho_0 \in \dot{B}_{p,1}^\frac{n}{p}.
$$
There exists a constant $c>0$ such that if $\|\rho_0-1\|_{\dot{B}_{p,1}^\frac{n}{p}}\leq c$ then there exists a positive time $T>0$ and a unique solution $(\rho, u, \nabla P)$ of \eqref{NSIK} with
$$
\begin{cases}
\rho-1, \mbox{ }\nabla \rho \in C_T(\dot{B}_{p,1}^\frac{n}{p}), \quad \d_t \rho \in L_T^\infty \dot{B}_{p,1}^{\frac{n}{p}-1} \cap L_T^2 \dot{B}_{p,1}^\frac{n}{p},\\
u\in C_T(\dot{B}_{p,1}^{\frac{n}{p}-1}), \quad (\d_t u, \nabla^2 u, \nabla P) \in L_t^1(\dot{B}_{p,1}^{\frac{n}{p}-1}).
\end{cases}
$$
Moreover, if in addition $\|u_0\|\leq c \mub$, then $T$ is bounded from below by some $T_{\kb}=\frac{C\mub}{\kb}$. Finally, if we denote by $(\rho_{\kb}, u_{\kb}, \n P_{\kb})$ the previous solution and by $(\rho, u, \n P)$  the global solution of the inhomogeneous incompressible Navier-Stokes system given by \cite{RP2} there exists a constant $C'>0$ such that if $T_{\kb}'=\frac{C\mub}{\kb^{1-\alpha}} < T_{\kb}$ for $0\leq \alpha \leq 1$:
$$
\begin{cases}
\vspace{0.2cm}
\|(u-u_{\kb}, \n (P- P_{\kb}))\|_{E_{T_{\kb}'}^{\frac{n}{p}-1}} \leq C' \kb^{\alpha},\\
\|\rho-\rho_{\kb}\|_{L_{T_{\kb}'}^\infty \dot{B}_{p,1}^\frac{n}{p}} \leq C' {\kb}^{\frac{3\alpha-1}{2}}, \mbox{ if } \frac{1}{3}\leq \alpha\leq 1.
\end{cases}
$$
where the space $E_T^s$ is defined in \eqref{espace}.
}
\label{T2}
\end{thm}

The article will be organized as follows. Sections 2 and 3 are devoted to the proof of the theorems. As the results adapt the methods from \cite{RP2, Dlagrangien2} we will skip details and mainly focus on what is new, namely the capillary term. The end of this section is devoted to some precisions and extensions of the results. Section 4 is a an appendix gathering definitions and main properties for homogeneous Besov spaces, Lagrangian change of variables, and estimates for the Lagrangian flow and its derivatives.

\section{Proof of Theorem \ref{T1}}

\subsection{Rescaling of the system}

If we introduce $(\trho, \tu, \tP)$ as follows:
\begin{equation}
\trho (t,x)=\rho(\frac{t}{\mub}, x),\quad \tu (t,x)=\frac{1}{\mub} u(\frac{t}{\mub}, x), \quad \tP (t,x)=\frac{1}{\mub^2} P(\frac{t}{\mub}, x).
\end{equation}
Then this allows us to study the case $(1,\kmu)$ instead of $(\mub, \kb)$ in the sense that $(\rho,u,\nabla P)$ solves \eqref{NSIK} if and only if $(\trho, \tu, \tP)$ solves:
\begin{equation}
\begin{cases}
\d_t \rho +u\cdot \n \rho=0,\\
\rho (\d_t u +u\cdot \n u) -\div \left(\mu (\rho) D(u)\right) +\n P= -\kmu \div \left(k(\rho) \n \rho \otimes \n \rho\right),\\
\div u=0,\\
(\rho,u)_{|t=0}=(\rho_0, \frac{1}{\mub} u_0).
\end{cases}
\label{NSIK2}
\tag{$NSIK2$}
\end{equation}
\begin{rem}
\sl{We emphasize that the ratio $\kmu$ also plays an important role in the compressible system, as observed for example in \cite{CH4}.}
\end{rem}
From now on we will focus on \eqref{NSIK2} and a solution for \eqref{NSIK2} on some $[0,T_1]$ will produce a solution for \eqref{NSIK} on $[0,\frac{T_1}{\mub}]$ thanks to the reverse change of variable:
\begin{equation}
\rho (t,x)=\trho (\mub t, x),\quad u (t,x)=\mub \tu(\mub t, x),\quad P (t,x)=\mub^2 \tP (\mub t, x).
\label{rescale}
\end{equation}

\subsection{Lagrangian formulation}

The method, first used in the incompressible case by T. Hmidi in \cite{TH1} then developped by Danchin (see for example \cite{Dlagrangien, Dlagrangien2}), Danchin and Mucha in \cite{RP2, RPbook}, is based on the following observation: if $(\rho, u, \nabla P)$ solves \eqref{NSIK2} and if $u$ is smooth enough (say Lipschitz), we introduce the flow $X$ associated to $u$, that is the solution to:
$$
\begin{cases}
\d_t X (t,x)=u(t,X(t,x)),\\
X(0,x)=x,
\end{cases}
$$
or equivalently $X(t,x)=x+\int_0^t u(\tau,X(\tau,x))d\tau$. As the jacobean determinant satisfies $\mbox{det}(DX(t,x))=e^{\int_0^t \div u(\tau, X(\tau,x)) d\tau}$, the incompressibility condition implies in our case $\mbox{det}(DX(t,.))\equiv 1$. Then, introducing for any function $f$, $\bar{f}(t,x)\overset{def}{=} f(t, X(t,x))$, with the computations from \cite{RP2} we obtain that $(\rhob, \ub, \n \Pb)$ solves the following system (we refer to the appendix for the transformation of the capillary term, the rest is dealt as in \cite{RP2, Dlagrangien2}):
\begin{equation}
\begin{cases}
\d_t \rhob=0,\\
\rhob \d_t \ub -\div \left[\mu (\rhob) A D_A (\ub)\right] +{^tA}\n \Pb= -\kmu \div \left[k(\rhob) A(^tA \n \rhob) \otimes (^tA \n \rhob)\right],\\
\div (A \ub)=0,\\
(\rhob,\ub)_{|t=0}=(\rho_0, \frac{1}{\mub} u_0),
\end{cases}
\end{equation}
where $D_A (z)=Dz\cdot A+^tA\cdot \n z$ and $A=DX^{-1}$. We emphasize that the new system presents two important features:
\begin{itemize}
\item there are no advection terms anymore,
\item the density becomes constant (this is due to incompressible setting, in the compressible setting the density is multiple of the inverse of the Jacobean determinant, we refer to \cite{Dlagrangien2} for more details).
\end{itemize}
The price to pay is reasonnable as it consists in dealing with Matrix $A$ which is close to $I_d$ when the time or the initial data are small as outlined in Proposition \ref{propestimL}. The previous system can then be rewritten into the form we will finally study (now denoting $(X_{\ub}, \Au)$ instead of $(X,A)$):
\begin{equation}
\begin{cases}
\rhob\equiv\ro,\\
\ro \d_t \ub -\div \left[\mu (\ro) \Au D_{\Au} (\ub)\right] +{^t}\Au \n \Pb= -\kmu \div \left[k(\ro) \Au (^t\Au \n \ro) \otimes (^t\Au \n \ro)\right],\\
\div (\Au \ub)=0,\\
\ub_{|t=0}=\frac{1}{\mub} u_0.
\end{cases}
\label{S1}
\end{equation}
\begin{rem}
\sl{Observe that with this notation $X_{\ub}(t,x)=x+\int_0^t \ub(\tau,x)d\tau$.}
\label{ChgtLagr}
\end{rem}
Next our stategy will be the same as in \cite{RP2, Dlagrangien2, CB}, let $\vb$ satisfying on $[0,T]$ for some $T>0$:
\begin{equation}
\begin{cases}
\vb \in \cC (\R_+, \dot{B}_{p,1}^{\frac{n}{p}-1}),\\
\d_t \vb, \n^2 \vb \in L_T^1 \dot{B}_{p,1}^{\frac{n}{p}-1},
\end{cases}
\mbox{ and} \quad
\int_0^T \|D \vb\|_{\dot{B}_{p,1}^{\frac{n}{p}-1}} \leq \ee \leq \ee_0,
\label{CondvT}
\end{equation}
where the small bound $\ee_0$ comes from \eqref{smallnes12} (see Proposition \ref{PropozitieFinal}). Let us define the flow:
$$
X_{\vb}(t,x)=x+\int_0^t \vb(\tau,x)d\tau,
$$
then all the work consists in proving that the following system (we recall that $adj(A)$ is the transposed cofactor matrix of $A$, that is $\mbox{det} A \cdot A^{-1}$ if $A$ is invertible):
\begin{equation}
\begin{cases}
\ro \d_t \ub -\div \left[\mu (\ro) \Av D_{\Av} (\ub)\right] +{^t}\Av \n \Pb= -\kmu \div \left[k(\ro) \Av (^t\Av \n \ro) \otimes (^t\Av \n \ro)\right],\\
\div (\adjv \ub)=0,\\
\ub_{|t=0}=\frac{1}{\mub} u_0,
\end{cases}
\label{Aux1}
\end{equation}
has a unique solution $(\ub, \n \Pb)$ on some $[0,T_0]$ for any $\vb$, which we reformulate in terms of fixed point: for any $(\vb, \n \Qb)$ the function $\Phi:(\vb, \n \Qb) \mapsto (\ub, \n \Pb)$ has a unique fixed point that solves \eqref{S1}. Next, thanks to the properties of $\ub$ we will be able to prove that $X_{\ub}$ is a $C^1$-diffeomorphism on $\R^n$ and we can perform the inverse change of variable $(\rho, u, P)=(\rhob, \ub, \Pb) \circ X_{\ub}^{-1}$. As a fixed point $\ub$ satisfies $\div (\adju \ub)$ which garantees that $\div u=0$ and as a consequence $X_{\ub}$ (which is the flow of $u$) is measure preserving and $\adju=\Au$, which implies that $\ub$ solves \eqref{S1}. Finally $(\rho, u, P)$ solves \eqref{NSIK2} on some $[0, T_0]$ which provides a unique solution of \eqref{NSIK} on the time interval $[0,\frac{T_0}{\mub}]$.
\begin{rem}
\sl{We emphasize that the given function $\vb$ genuinely depends on the Lagrangian variables, it is not the change of variable of a given function $v$. In other words, the system is solved in Lagrangian variables then the inverse change of variables is performed.}
\end{rem}

\subsection{A priori estimates}

As in \cite{RP2} the main ingredient in the proof of Theorem \ref{T2} is the following maximal regularity estimates proved in \cite{RP2, RPbook}:

\begin{prop}
\sl{Let $p\in[1,\infty]$, $s\in \R$, $u_0 \in \dot{B}_{p,1}^s(\R^n)$ and $f\in L_T^1 \dot{B}_{p,1}^s(\R^n)$. Let $M:[0,T]\times \R^n \rightarrow \R$ such that:
$$
\begin{cases}
\n \div M, \mbox{ }\d_t M \in L_T^1 \dot{B}_{p,1}^s(\R^n),\\
\mathcal{Q} M \in \cC([0,\infty), \dot{B}_{p,1}^{\frac{n}{p}-1}(\R^n)),\\
\div M|_{t=0}=\div u_0 \mbox{ on }\R^n.
\end{cases}
$$
Then the system
\begin{equation}
\begin{cases}
\d_t u-\mu \Delta u +\n P=f,\\
\div u=\div M,\\
u|_{t=0}=u_0,
\end{cases}
\end{equation}
admits a unique solution $(u,\n P)$ in the following space
$$
E_T^s=\left\{(u,\n P), \mbox{ } u\in \cC_T \dot{B}_{p,1}^s(\R^n) \mbox{ and } \d_t u, \n^2 u, \n P \in L_T^1 \dot{B}_{p,1}^s(\R^n)\right\},
$$
and there exists a positive constant $C$ (independant of $\mu, T$) such that
\begin{multline}
\|(u,\n P)\|_{E_T^s} \overset{def}{=} \|u\|_{L_T^\infty \dot{B}_{p,1}^s}+ \|(\d_t u, \mu \n^2 u, \n P)\|_{L_T^1 \dot{B}_{p,1}^s}\\
\leq C \left( \|u_0\|_{\dot{B}_{p,1}^s}+\|(f, \mu \n \div M, \d_t M)\|_{L_T^1 \dot{B}_{p,1}^s}\right)
\label{espace}
\end{multline}
\label{Apriori1}
}
\end{prop}
Similarly, the main ingredient in the proof of Theorem \ref{T1} (no smallness assumptions or additional low frequency regularity) is the new estimates recently obtained by the first author in \cite{CB}:
\begin{prop}
\sl{Let $n=2$ and $p\in(1,4)$ or $n=3$ and $p\in(\frac{6}{5},4)$. Let $a,b$ two functions such that there exists positive constants $(a_*, a^*, b_*, b^*, \bar{a}, \bar{b})$ such that
$$
\begin{cases}
a-\bar{a}, b-\bar{b} \in \dot{B}_{p,1}^\frac{n}{p}(\R^n),\\
0<a_* \leq a \leq a^*,\\
0<b_* \leq b \leq b^*.
\end{cases}
$$
Let $u_0 \in \dot{B}_{p,1}^{\frac{n}{p}-1}(\R^n)$ and $f\in L_T^1 \dot{B}_{p,1}^{\frac{n}{p}-1}(\R^n)$. Let $M:[0,T]\times \R^n \rightarrow \R$ such that:
$$
\begin{cases}
\n \div M, \mbox{ }\d_t M \in L_T^1 \dot{B}_{p,1}^{\frac{n}{p}-1}(\R^n),\\
\mathcal{Q} M \in \cC([0,\infty), \dot{B}_{p,1}^{\frac{n}{p}-1}(\R^n)),\\
\div M|_{t=0}=\div u_0 \mbox{ on }\R^n.
\end{cases}
$$
Then System
\begin{equation}
\begin{cases}
\d_t u-a \div(bD(u)) +a\n P=f,\\
\div u=\div M,\\
u|_{t=0}=u_0,
\end{cases}
\end{equation}
admits a unique solution $(u,\n P)\in E_T=E_T^{\frac{n}{p}-1}$ and there exists a positive constant $C=C(a,b)$ such that for all $t\leq T$,
$$
\|(u,\n P)\|_{E_t} \leq e^{C(t+1)}\left( \|u_0\|_{\dot{B}_{p,1}^{\frac{n}{p}-1}}+\|(f, \mu \n \div M, \d_t M)\|_{L_t^1 \dot{B}_{p,1}^{\frac{n}{p}-1}}\right).
$$
\label{Apriori2}
}
\end{prop}

\subsection{First step: well-posedness for \eqref{Aux1}}
With a aim of conciseness, we will only give details for terms and computations related to the capillary term as the other terms are the same as in \cite{CB, RP2}.
\\
Let $\vb$ and $T$ given as in \eqref{CondvT}, for $(\wb,\n \Qb)\in E_T$ let us consider the following system:
\begin{equation}
\begin{cases}
\d_t \ub -\frac{1}{\ro}\div \left[\mu (\ro)D(\ub)\right] +\frac{1}{\ro} \n \Pb= F_{ext} =F_{\vb}^1(\n \Qb) +F_{\vb}^2(\wb) +F_{\vb}^3 ,\\
\div \ub= \div M_{\vb}(\wb),\\
\ub_{|t=0}=\frac{1}{\mub} u_0,
\end{cases}
\label{Aux2}
\end{equation}
where
$$
\begin{cases}
F_{\vb}^1(\n \Qb) \overset{def}{=} \frac{1}{\ro}(I_d- ^t\Av) \n \Qb,\\
F_{\vb}^2(\wb)\overset{def}{=} \frac{1}{\ro} \div \Big[\mu (\ro)(\Av D_{\Av}(\wb)-D(\wb))\Big],\\
F_{\vb}^3\overset{def}{=} -\kmu \frac{1}{\ro} \div \Big[k(\ro) \Av (^t\Av \n \ro) \otimes (^t\Av \n \ro)\Big],\\
M_{\vb}(\wb) \overset{def}{=} (I_d-\adjv)\wb.
\end{cases}
$$
Proposition \ref{Apriori2} implies there exists a unique solution $(\ub,\n \Pb)\in E_T$ and a constant $C_{\ro}$ such that for all $t\leq T$,
\begin{multline}
\|(\ub,\n \Pb)\|_{E_t} =\|u\|_{L_t^\infty \dot{B}_{p,1}^{\frac{n}{p}-1}}+ \|(\d_t u, \n^2 u, \n P)\|_{L_t^1 \dot{B}_{p,1}^{\frac{n}{p}-1}}\\
\leq e^{C_{\ro}(t+1)}\left( \frac{1}{\mub}\|u_0\|_{\dot{B}_{p,1}^{\frac{n}{p}-1}}+\|(F_{ext}, \n \div M_{\vb}(\wb), \d_t M_{\vb}(\wb))\|_{L_t^1 \dot{B}_{p,1}^{\frac{n}{p}-1}}\right)
\end{multline}
The first two external terms are the same as in \cite{RP2, CB}, and using Condition \eqref{CondvT} and \eqref{A1} from Proposition \ref{propestimL}, we easily obtain that:
$$
\begin{cases}
\|F_{\vb}^1(\n \Qb)\|_{L_t^1 \dot{B}_{p,1}^{\frac{n}{p}-1}} \leq C_{\ro} \|D\vb\|_{L_t^1 \dot{B}_{p,1}^\frac{n}{p}} \|\n \Qb\|_{L_t^1 \dot{B}_{p,1}^{\frac{n}{p}-1}},\\
\|F_{\vb}^2(\wb)\|_{L_t^1 \dot{B}_{p,1}^{\frac{n}{p}-1}} \leq C_{\ro} (1+\|D\vb\|_{L_t^1 \dot{B}_{p,1}^\frac{n}{p}}) \|D\vb\|_{L_t^1 \dot{B}_{p,1}^\frac{n}{p}} \|D\wb\|_{L_t^1 \dot{B}_{p,1}^\frac{n}{p}},
\end{cases}
$$
The third term is dealt with the classical product laws in Besov spaces (see section \ref{appendiceLP}), using the fact that $k(\ro)=1+k(\ro)-k(1)$ (recall that $k(1)=1$) and \eqref{A1}:
\begin{multline}
\|F_{\vb}^3\|_{L_t^1 \dot{B}_{p,1}^{\frac{n}{p}-1}} \leq C \kmu (1+\|\frac{1}{\ro}-1\|_{\dot{B}_{p,1}^\frac{n}{p}}) (1+\|k(\ro)-1\|_{\dot{B}_{p,1}^\frac{n}{p}}) \|\Av (^t\Av \n \ro) \otimes (^t\Av \n \ro)\|_{L_t^1 \dot{B}_{p,1}^{\frac{n}{p}-1}}\\
\leq C_{\ro} \kmu t(1+\|\Av-I_d\|_{L_t^\infty \dot{B}_{p,1}^\frac{n}{p}}) \left((1+\|(^t \Av-I_d)\|_{L_t^\infty \dot{B}_{p,1}^\frac{n}{p}}) \| \n \ro\|_{\dot{B}_{p,1}^\frac{n}{p}}\right)^2\\
\leq C_{\ro} \kmu t (1+\|D\vb\|_{L_t^1 \dot{B}_{p,1}^\frac{n}{p}}) ^3.
\end{multline}
Thanks to \eqref{magic}, we can write
$$
\div M_{\vb}(\wb) = \div [(I_d-\adjv)\wb]=D\wb:(I_d-\Jv \Av),
$$
and thanks to Proposition \ref{propestimL} we obtain as in \cite{CB}:
$$
\begin{cases}
\|\n \div M_{\vb}(\wb)\|_{L_t^1 \dot{B}_{p,1}^{\frac{n}{p}-1}} \leq C \|D\vb\|_{L_t^1 \dot{B}_{p,1}^\frac{n}{p}} \|D\wb\|_{L_t^1 \dot{B}_{p,1}^\frac{n}{p}},\\
\|\d_t M_{\vb}(\wb)\|_{L_t^1 \dot{B}_{p,1}^{\frac{n}{p}-1}} \leq C \|D\vb\|_{L_t^1 \dot{B}_{p,1}^\frac{n}{p}} (\|D\wb\|_{L_t^\infty \dot{B}_{p,1}^{\frac{n}{p}-1}} +\|\d_t\wb\|_{L_t^1 \dot{B}_{p,1}^{\frac{n}{p}-1}}).
\end{cases}
$$
So that according to condition \eqref{CondvT} we end up for all $t\leq T$ with,
$$
\|(\ub,\n \Pb)\|_{E_t} \leq e^{C_{\ro}(t+1)} \left( \frac{1}{\mub}\|u_0\|_{\dot{B}_{p,1}^{\frac{n}{p}-1}} +C_{\ro} \ee \|(\wb, \n \Qb)\|_{E_t} +C_{\ro} \kmu t \right),
$$
and the application $\Psi: (\wb, \n \Qb) \mapsto (\ub,\n \Pb)$ is well defined $E_T \rightarrow E_T$. Next let us prove that this is a contraction: given $(\wb_i, \n \Qb_i)\in E_T$ ($i=1,2$), if $(\ub_i,\n \Pb_i) = \Psi (\vb_i,\n \Qb_i)$, then $(\dub, \n \dPb) \overset{def}{=} (\ub_2-\ub_1, \n \Pb_2- \n \Pb_1)$ solves:
\begin{equation}
\begin{cases}
\d_t \dub -\frac{1}{\ro}\div \left[\mu (\ro)D(\dub)\right] +\frac{1}{\ro} \n \dPb =F_{\vb}^1(\n \dQb) +F_{\vb}^2(\dwb) ,\\
\div \dub= \div M_{\vb}(\dwb),\\
\dub_{|t=0}=0.
\end{cases}
\end{equation}
Similar computations as before imply that for all $t\leq T$,
$$
\|(\dub,\n \dPb)\|_{E_t} \leq C_{\ro} \ee e^{C_{\ro}(t+1)} \|(\dwb, \n \dQb)\|_{E_t},
$$
and $\Psi$ is $\frac{1}{2}$-Lipschitz when for example $T\leq 1$ and $\ee<\frac{1}{2C_{\ro}} e^{-2C_{\ro}}$ and has a unique fixed point $(\ub,\n \Pb) \in E_T$ which solves \eqref{Aux1}. Moreover the fixed point belongs to $E_T$ and for all $t\leq T$,
$$
\|(\ub,\n \Pb)\|_{E_t} \leq 2e^{C_{\ro}(t+1)} \left( \frac{1}{\mub}\|u_0\|_{\dot{B}_{p,1}^{\frac{n}{p}-1}} +C_{\ro} \kmu t \right),
$$
\begin{rem}
\sl{We could also ask $T\leq B$ and $\ee<\frac{1}{2C_{\ro}} e^{-(1+B)C_{\ro}}$ for some $B>0$, but as we will see later, taking a large $B$ will not help for the final lifespan.
}
\label{remB}
\end{rem}
Let us introduce the space
$$
F_T^\ee=\{(\vb, \n \Qb)\in E_T \mbox{ with } \|\n \vb\|_{L_T^1 \dot{B}_{p,1}^\frac{n}{p}} \leq \ee\}.
$$
We just proved that given any $(\vb, \n \Qb)\in F_T^\ee$ and under Condition \eqref{CondvT}, System \eqref{Aux1} has a unique solution $(\ub, \n \Pb)=\Phi(\vb, \n \Qb)$. Let us recall that our aim is to prove that \eqref{S1} has a unique solution by proving that $\Phi$ has a unique fixed point. Proving that $\Phi$ is a contraction will not be difficult but the problem is that as $u_0$ is not assumed to be small, there is no reason for $\ub$ to satisfy the following condition (which is crucial to go back to the original variables):
$$
\|\n \ub\|_{L_T^1 \dot{B}_{p,1}^\frac{n}{p}} \leq \ee \leq \ee_0,
$$
even if there exists some $T_{\ub}>0$ such that the integral is bounded by $\ee$, we do not know if $T_{\ub} \geq T$. In other words, we are not able to prove that $\Phi$ maps $F_T^\ee$ into itself due to this integral condition.

\subsection{Second step: well-posedness for \eqref{S1}}

As in \cite{RP2, CB}, to overcome this problem, the idea is to introduce the free solution, let us define $(u_L, \n P_L)$ the unique global solution of the following system:
\begin{equation}
\begin{cases}
\d_t u_L -\frac{1}{\ro}\div \left[\mu (\ro)D(u_L)\right] +\frac{1}{\ro} \n P_L =-\kmu \frac{1}{\ro} \div\left[ k(\ro) \n \ro \otimes \n \ro\right] \overset{def}{=}F_0^3 ,\\
\div u_L =0,\\
\ub_{|t=0}=\frac{1}{\mub} u_0.
\end{cases}
\label{SystL}
\end{equation}
The a priori estimates also provide the fact that for all $t$,
$$
\|(u_L,\n P_L)\|_{E_t} \leq e^{C_{\ro}(t+1)} \left( \frac{1}{\mub}\|u_0\|_{\dot{B}_{p,1}^{\frac{n}{p}-1}} +C_{\ro} \kmu t \right).
$$
We are now able to precise the parameters: assume that $\ee$ and $T$ satisfy:
\begin{equation}
\begin{cases}
\vspace{0.2cm}
\ee= \min(\ee_0, \frac{e^{-2C_{\ro}}}{32 C_{\ro}}),\\
T=\min \left(1, \frac{\mub^2}{\kb} \frac{ e^{-2C_{\ro}}}{32 C_{\ro}}, \sup \Big\{ t>0, \mbox{ } \|(\d_t u_L, \n^2 u_L, \n P_L)\|_{L_t^1 \dot{B}_{p,1}^{\frac{n}{p}-1}} +\|u_L\|_{L_t^2 \dot{B}_{p,1}^\frac{n}{p}} \leq \frac{\ee}{2}\Big\}\right),
\end{cases}
\label{CondepsT}
\end{equation}
let us define the space
$$
G_T^\ee=\Big\{(f, \n g)\in E_T \mbox{ with } \|(f,\n g)\|_{E_T}\leq \frac{\ee}{2}\Big\}.
$$
Thanks to Condition \eqref{CondepsT}, for any $(\tv,\n \tQ)\in G_T^\ee$, then $(\vb, \n \Qb)\overset{def}{=} (u_L+\tv, \n P_L + \n \tP)$ belongs to $F_T^\ee$ and we easily check that Condition \eqref{CondvT} is satisfied so there exists a unique solution $(\ub, \n \Pb)=\Phi(\vb, \n \Qb)$. If we set $(\tu, \n \tP)=(\ub-u_L, \n \Pb -\n P_L)$ then if $T,\ee$ are small enough, we can prove that $(\tu, \n \tP) \in G_T^\ee$: indeed they satisfy the following system
\begin{equation}
\begin{cases}
\d_t \tu -\frac{1}{\ro}\div \left[\mu (\ro)D(\tu)\right] +\frac{1}{\ro} \n \tP= F_{ext}' =F_{\vb}^1(\n \Pb) +F_{\vb}^2(\ub) +G_{\vb}^3 ,\\
\div \tu= \div M_{\vb}(\ub),\\
\tu_{|t=0}=0,
\end{cases}
\label{Aux3}
\end{equation}
where $A_{\vb}^1(\n \Pb)$, $A_{\vb}^2(\ub)$ and $M_{\vb}(\ub)$ are the same as in \eqref{Aux2}, and the last term is:
\begin{multline}
G_{\vb}^3 \overset{def}{=} F_{\vb}^3-F_0^3=-\kmu \frac{1}{\ro} \div \Big[k(\ro) \Big\{\Av \Big((^t\Av \n \ro) \otimes \big((^t\Av -I_d) \n \ro\big)\\
+\big((^t\Av -I_d) \n \ro\big) \otimes \n \ro\Big) +(\Av-I_d) \n \ro \otimes \n \ro \Big\}\Big].
\end{multline}
As before, we obtain thanks to Propositions \ref{Apriori2} and \ref{propestimL}:
\begin{equation}
\|F_{\vb}^1(\n \Pb)\|_{L_t^1 \dot{B}_{p,1}^{\frac{n}{p}-1}} \leq C_{\ro} (\|D\tv\|_{L_t^1 \dot{B}_{p,1}^\frac{n}{p}} +\|D u_L\|_{L_t^1 \dot{B}_{p,1}^\frac{n}{p}}) (\|\n \tP\|_{L_t^1 \dot{B}_{p,1}^{\frac{n}{p}-1}} +\|\n P_L\|_{L_t^1 \dot{B}_{p,1}^{\frac{n}{p}-1}} ),
\end{equation}
\begin{multline}
\|F_{\vb}^2(\wb)\|_{L_t^1 \dot{B}_{p,1}^{\frac{n}{p}-1}} \leq C_{\ro} (1+\|D\tv\|_{L_t^1 \dot{B}_{p,1}^\frac{n}{p}} +\|D u_L\|_{L_t^1 \dot{B}_{p,1}^\frac{n}{p}})\\
\times (\|D\tv\|_{L_t^1 \dot{B}_{p,1}^\frac{n}{p}} +\|D u_L\|_{L_t^1 \dot{B}_{p,1}^\frac{n}{p}}) (\|D\tu\|_{L_t^1 \dot{B}_{p,1}^\frac{n}{p}} +\|D u_L\|_{L_t^1 \dot{B}_{p,1}^\frac{n}{p}}) ,
\end{multline}
\begin{multline}
\|G_{\vb}^3\|_{L_t^1 \dot{B}_{p,1}^{\frac{n}{p}-1}} \leq C_{\ro} \kmu t \|\n \ro\|_{\dot{B}_{p,1}^\frac{n}{p}}^2 (1+\|D\tv\|_{L_t^1 \dot{B}_{p,1}^\frac{n}{p}} +\|D u_L\|_{L_t^1 \dot{B}_{p,1}^\frac{n}{p}})^2\\
\times (\|D\tv\|_{L_t^1 \dot{B}_{p,1}^\frac{n}{p}} +\|D u_L\|_{L_t^1 \dot{B}_{p,1}^\frac{n}{p}}),
\end{multline}
and using once again $\div M_{\vb}(\ub) = \div [(I_d-\adjv)\ub]=D\ub:(I_d-\Jv \Av)$ and Proposition \ref{propestimL}, we get:
$$
\|\n \div M_{\vb}(\wb)\|_{L_t^1 \dot{B}_{p,1}^{\frac{n}{p}-1}} \leq C (\|D\tu\|_{L_t^1 \dot{B}_{p,1}^\frac{n}{p}} +\|D u_L\|_{L_t^1 \dot{B}_{p,1}^\frac{n}{p}}) (\|D\tv\|_{L_t^1 \dot{B}_{p,1}^\frac{n}{p}} +\|D u_L\|_{L_t^1 \dot{B}_{p,1}^\frac{n}{p}}).
$$
As we wish to take advantage of the smallness of the $L_T^1$ norm, we cannot do as before for the last term (as the $\|u_L\|_{L_T^\infty \dot{B}_{p,1}^{\frac{n}{p}-1}}$ norm may not be small and would prevent any absorption by the left-hand side) and the idea is (as in \cite{RP2, CB}) to use the $L_t^2 \dot{B}_{p,1}^\frac{n}{p}$-norm of $u_L$ instead:
\begin{multline}
\|\d_t M_{\ub}(\wb)\|_{L_t^1 \dot{B}_{p,1}^{\frac{n}{p}-1}} \leq C \|\d_t \adjv\|_{L_T^2 \dot{B}_{p,1}^{\frac{n}{p}-1}} \|\ub\|_{L_T^2 \dot{B}_{p,1}^\frac{n}{p}} +\|I_d -\adjv\|_{L_t^\infty \dot{B}_{p,1}^\frac{n}{p}} \|\d_t \ub\|_{L_t^1 \dot{B}_{p,1}^{\frac{n}{p}-1}}\\
\leq C (\|D\tv\|_{L_t^2 \dot{B}_{p,1}^{\frac{n}{p}-1}} +\|D u_L\|_{L_t^2 \dot{B}_{p,1}^{\frac{n}{p}-1}}) (\|\tu\|_{L_t^2 \dot{B}_{p,1}^\frac{n}{p}} +\|u_L\|_{L_t^2 \dot{B}_{p,1}^\frac{n}{p}})\\
+C (\|D\tv\|_{L_t^1 \dot{B}_{p,1}^\frac{n}{p}} +\|D u_L\|_{L_t^1 \dot{B}_{p,1}^\frac{n}{p}}) (\|\d_t \tu\|_{L_t^1 \dot{B}_{p,1}^{\frac{n}{p}-1}} +\|\d_t u_L\|_{L_t^1 \dot{B}_{p,1}^{\frac{n}{p}-1}}).
\end{multline}
Thanks to the fact that $(\tv, \n \tQ)\in G_T^\ee$ and condition \eqref{CondepsT} we end up for all $t\leq T$ with (we do not give details as the computations are the same as before),
$$
\|(\tu,\n \tP)\|_{E_t} \leq \ee C_{\ro} e^{C_{\ro}(t+1)} \left(\|(\tu,\n \tP)\|_{E_t} +\ee +\kmu T \right),
$$
and then, as $\ee C_{\ro} e^{C_{\ro}(T+1)} \leq \frac{1}{2}$ (thanks to \eqref{CondepsT}),
$$
\|(\tu,\n \tP)\|_{E_t} \leq 2 \ee C_{\ro} e^{C_{\ro}(T+1)} \left(\ee +\kmu T \right) \leq \frac{\ee}{2}.
$$
Finally, we proved that for any $(\tv,\n \tQ)\in G_T^\ee$ there exists a unique solution $(\tu,\n \tP) \in G_T^\ee$ of System \eqref{Aux3}. We will abusively denote once more this solution by $\Phi (\tv,\n \tQ)$. To prove that $\Phi$ is a contraction is very close to what we did for the application $\Psi$ so we will not give details. All that remains is to get back to the original variables: as the unique fixed point of $\Phi$, $\ub$ satisfies $\div(\adju \ub)=0$ and as we also have $\|\n \ub\|_{L_T^1 \dot{B}_{p,1}^\frac{n}{p}} \leq \ee \leq \ee_0$, Proposition \ref{PropozitieFinal} implies that $X_{\ub}$ is a global measure-preserving $C^1$-diffeomorphism on $\R^n$. As a consequence of Proposition \ref{Comp}, $(\rho, u, P)=(\bar{\rho}, \ub, \Pb) \circ (X_{\ub})^{-1}$ has the announced regularities and solves \eqref{NSIK2}: for example as $\ro-1$ and $\n \ro \in \dot{B}_{p,1}^\frac{n}{p}$, so do $\rho-1$ and $\n \rho$ and as $\d_t \rho =-u\cdot \n \rho$ we obtain the $L_T^\infty \dot{B}_{p,1}^{\frac{n}{p}-1}$ and $L_T^2 \dot{B}_{p,1}^\frac{n}{p}$ estimates for $\d_t \rho$.

Uniqueness for System \eqref{S1} concludes the proof of Theorem \ref{T1}. $\blacksquare$
\\

Let us end this section with precisions about the dependancy of $T$ in terms of the physical parameters: a more precise use of Proposition \ref{Apriori2} gives that $(u_L, \n P_L)$, solution of \eqref{SystL}, satisfies:
$$
\|(u_L,\n P_L)\|_{E_t} \leq 2e^{C_{\ro}(t+1)} \left( \frac{1}{\mub}\Sum_{j\in \Z} (1-e^{-CT2^{2j}}) 2^{j(\frac{n}{p}-1)} \|\ddj u_0\|_{L^p} +C_{\ro} \kmu t \right),
$$
and due to the exponential term we can see that if we asked $T\leq B$ for some large $B$ as in Remark \ref{remB}, the above estimates implies that $t\leq \frac{\mub^2}{\kb} \frac{\ee e^{-(1+B) C_{\ro}}}{4C_{\ro}}$, which explains why $B=1$ is sufficient as outlined in Remark \ref{remB}.

\section{Proof of Theorem \ref{T2}}

Let us now turn to the cases when $\ro-1$ is small (allowing the full ranges for $n,p$). As announced in the previous section, instead of System \eqref{Aux1}, we will take advantage of Proposition \ref{Apriori1} (As $\div D(u)=\Delta u +\n \div u$ we are in the scope of this result) and study:
\begin{equation}
\begin{cases}
\d_t \ub -\div \left[D(\ub)\right] + \n \Pb= H^1(\ub)+H_{\vb}^2(\n \Pb) +H_{\vb}^3(\ub) +H_{\vb}^4,\\
\div \ub= \div M_{\vb}(\ub),\\
\ub_{|t=0}=\frac{1}{\mub} u_0,
\end{cases}
\label{Aux4}
\end{equation}
where $M_{\vb}(\ub)$ is the same as before and the external force terms are very close to those from \eqref{Aux2} (recall that $\mu(1)=1$):
$$
\begin{cases}
H^1(\wb) \overset{def}{=} (1-\ro) \d_t \wb,\\
H_{\vb}^2(\n \Qb) \overset{def}{=} (I_d- ^t\Av) \n \Qb,\\
H_{\vb}^3(\ub) \overset{def}{=} \div \left[\mu (\ro)\Av D_{\Av}(\wb)-\mu(1)D(\wb)\right],\\
H_{\vb}^3\overset{def}{=} -\kmu \div \left[k(\ro) \Av (^t\Av \n \ro) \otimes (^t\Av \n \ro)\right],\\
\end{cases}
$$
As before, if $(\vb, T)$ satisfies Condition \eqref{CondvT}, if $(\wb, \n \Qb)\in E_T$, let us introduce the following system:
\begin{equation}
\begin{cases}
\d_t \ub -\div \left[D(\ub)\right] + \n \Pb= H^1(\wb)+H_{\vb}^2(\n \Qb) +H_{\vb}^3(\wb) +H_{\vb}^4,\\
\div \ub= \div M_{\vb}(\wb),\\
\ub_{|t=0}=\frac{1}{\mub} u_0.
\end{cases}
\label{Aux5}
\end{equation}
Thanks to Proposition \ref{Apriori1}, with similar computations as before, we obtain that for all $t\leq T$,
$$
\|(\ub,\n \Pb)\|_{E_t} \leq C \left( \frac{1}{\mub}\|u_0\|_{\dot{B}_{p,1}^{\frac{n}{p}-1}} +\big(\|1-\ro\|_{\dot{B}_{p,1}^\frac{n}{p}} +\ee\big) \|(\wb, \n \Qb)\|_{E_t} +C_{\ro} \kmu t \right),
$$
and the application $\Psi: (\wb, \n \Qb) \mapsto (\ub,\n \Pb)$ is well defined $E_T \rightarrow E_T$. Similarly as before we easily show this is a contraction if $C(\|1-\ro\|_{\dot{B}_{p,1}^\frac{n}{p}} +\ee) \leq \frac{1}{2}$ for example and then we obtain a unique fixed point satisfying for all $t\leq T$:
$$
\|(\ub,\n \Pb)\|_{E_t} \leq 2C \left( \frac{1}{\mub}\|u_0\|_{\dot{B}_{p,1}^{\frac{n}{p}-1}} +C_{\ro} \kmu t \right),
$$
\begin{rem}
\sl{Contrary to the result from \cite{RP2}, a smallness condition only on $\|\ro-1\|_{\mathcal{M}(\dot{B}_{p,1}^{\frac{n}{p}-1})}$ will not be sufficient because the capillary term features $(k(\ro)-k(1)$ to be estimated in Besov spaces.}
\end{rem}
We then have to split in two cases wether $u_0$ is small or not.

\textbf{First case :} if $u_0$ and $T$ are so small that
$$
2C \left( \frac{1}{\mub}\|u_0\|_{\dot{B}_{p,1}^{\frac{n}{p}-1}} +C_{\ro} \kmu T \right) \leq \ee,
$$
that is for example when
\begin{equation}
\|u_0\|_{\dot{B}_{p,1}^{\frac{n}{p}-1}} \leq \frac{\ee}{4C} \mub \quad \mbox{and } T\leq \frac{\ee}{4CC_{\ro}} \frac{\mub^2}{\kb},
\label{CondT0}
\end{equation}
then $\|(\ub,\n \Pb)\|_{E_T} \leq \ee$ and the application $\Phi$, associating to $(\vb, \n \Qb)$ the unique solution $(\ub, \n \Pb)$ is well defined $F_T^\ee \rightarrow F_T^\ee$. As no free system is required (and no time exponential appears in the a priori estimates) we can simply take $T=\frac{\ee}{4CC_{\ro}} \frac{\mub^2}{\kb}$. We prove similarly that $\Phi$ is a contraction and then we obtain a unique fixed point for $\Phi$ which provides the unique solution for System \eqref{S1} and concludes, thanks to the inverse Lagrangian change of variable and inverse scale change, the proof of the last part of Theorem \ref{T2} (giving a lifespan bounded from below by a multiple of $\frac{\mub}{\kb}$). $\blacksquare$
\\

\textbf{Second case :} if $u_0$ is not assumed to be small, due to the $L_T^\infty$-norm nothing garantees that $(\ub,\n \Pb)\in F_T^\ee$ anymore and the idea is, as before, to introduce the unique solution $(u_L, \n P_L)$ of the free system (simpler than the one in the previous section):
\begin{equation}
\begin{cases}
\d_t u_L -\div \left[D(u_L)\right] +\n P_L =-\kmu \div\left[ k(\ro) \n \ro \otimes \n \ro\right] =H_0^4,\\
\div u_L =0,\\
\ub_{|t=0}=\frac{1}{\mub} u_0.
\end{cases}
\label{SystL2}
\end{equation}
With the same arguments as in the previous section we prove the existence of a constant $C_{\ro}>0$ such that if:
\begin{equation}
\begin{cases}
\vspace{0.2cm}
\ee= \min(\ee_0, \frac{1}{32 C_{\ro}}),\\
T=\min \left(1, \frac{\mub^2}{\kb} \frac{1}{32 C_{\ro}}, \sup \Big\{ t>0, \mbox{ } \|(\d_t u_L, \n^2 u_L, \n P_L)\|_{L_t^1 \dot{B}_{p,1}^{\frac{n}{p}-1}} +\|u_L\|_{L_t^2 \dot{B}_{p,1}^\frac{n}{p}} \leq \frac{\ee}{2}\Big\}\right),
\end{cases}
\label{CondepsT2}
\end{equation}
then (similarly to what we did in the previous section, we leave the details to the reader), the application $\Phi$ mapping $(\tv, \n \tQ)$ to the unique solution $(\tu, \n \tP)$ of the following system:
\begin{equation}
\begin{cases}
\d_t \tu -\div \left[D(\tu)\right] + \n \tP= H^1(u_L+\tu)+H_{u_L+\tv}^2(\n P_L+\n \tQ) +H_{u_L+\tv}^3(u_L+\tu)\\
\hspace{5cm}
 +(H_{u_L+\tv}^4 -H_0^3),\\
\div \tu= \div M_{u_L+\tv}(u_L+\tu),\\
\tu_{|t=0}=0,
\end{cases}
\label{Aux6}
\end{equation}
is well defined from $G_T^\ee$ to itself, and contractive. Then there exists a unique fixed point which ends the proof of the rest of Theorem \ref{T2}. $\blacksquare$

\subsection{Convergence when $\kb$ goes to zero}

If $(\rho, u, P)$ and $(\rhok, u_{\kb}, P_{\kb})$ are defined as in Theorem \ref{T2}, performing both of the Lagrangian changes of variables, we obtain that the difference $(\du, \n \dP) \overset{def}{=} (\uk-\ub, \Pk-\Pb)$ satisfies the following system:
\begin{equation}
\begin{cases}
\d_t \du -\div \left[D(\du)\right] + \n \Pk= H^1(\du)+H_{\ub}^2(\n \dP) +\Sum_{i=1}^3 K^i +H_{\uk}^3,\\
\div \du= \div M,\\
{\du}_{|t=0}=0,
\end{cases}
\end{equation}
where all the right-hand side terms are the same as in \eqref{Aux4} except:
$$
\begin{cases}
K^1 \overset{def}{=} -(^tA_{\uk}-^tA_{\ub})\cdot \n \Pk,\\
K^2 \overset{def}{=} \div \left[\mu(\ro)\Big((A_{\uk}-A_{\ub}) \cdot D_{A_{\uk}}(\uk)-D_{A_{\uk}-A_{\ub}}(\ub) \Big)\right],\\
K^3 \overset{def}{=} \div \left[\mu(\ro)(A_{\ub} \cdot D_{A_{\uk}}(\du) -D \du) +(\mu(\ro)-1) D\du\right],\\
\div M= \div \Big[-(A_{\uk}-A_{\ub}) \uk +(I_d-A_{\ub}) \du\Big] =-D\uk :(A_{\uk}-A_{\ub})+ D\du: (I_d-A_{\ub}).
\end{cases}
$$
Using the same arguments as before we obtain that under Condition \eqref{CondT0}, for all $t\leq T_{\kb}\overset{def}{=} \frac{\ee \mub^2}{4C C_{\ro}} \frac{1}{\kb}$,
$$
\|(\du, \n \dP)\|_{E_t} \leq \ee \|(\du, \n \dP)\|_{E_t} + C_{\ro} \kb t,
$$
so that (as $\ee \leq \frac{1}{2}$) for all $t\leq T_{\kb}'\overset{def}{=} \frac{\ee \mub^2}{4C C_{\ro}} \frac{1}{\kb^{1-\aa}}$ we end up with:
$$\|(\du, \n \dP)\|_{E_{T_{\kb}'}} \leq C \kb^{\alpha}.$$
And for the density as for all $t,x$ :
$$
\delta \rho(t,x) =\rho_{\kb}(t,x)-\rho(t,x)=\ro (X_{\uk}(t,x))-\ro (X_{\ub}(t,x)),
$$
then
$$
\|\delta \rho\|_{L_{T_{\kb}'}^\infty \dot{B}_{p,1}^\frac{n}{p}} \leq \|\n \ro\|_{\dot{B}_{p,1}^\frac{n}{p}} \|\du\|_{L_{T_{\kb}'}^1 \dot{B}_{p,1}^\frac{n}{p}} \leq C_{\ro} \sqrt{T_{\kb}'} \|\du\|_{L_t^2 \dot{B}_{p,1}^\frac{n}{p}} \leq C_{\ro} \kb (T_{\kb}')^\frac{3}{2} \leq C_{\ro} \kb^{\frac{3\alpha -1}{2}},
$$
which ends the proof of the theorem. $\blacksquare$

\subsection{Precisions about the lifespan}

In both proofs, we had to introduce, for some fixed small $\ee$ (smaller than $\ee_0$ from \eqref{smallnes12} from Proposition \ref{PropozitieFinal}), the time:
\begin{equation}
T_\ee \overset{def}{=} \sup \{ t>0, \mbox{ } \|(\d_t u_L, \n^2 u_L, \n P_L)\|_{L_t^1 \dot{B}_{p,1}^{\frac{n}{p}-1}} +\|u_L\|_{L_t^2 \dot{B}_{p,1}^\frac{n}{p}} \leq \frac{\ee}{2}\}.
\label{CondT2}
\end{equation}
Let us give more details for example in the second case (small $\rho_0-1$): as $(u_L, \n P_L)$ solves System \eqref{SystL2}, projecting thanks to the Leray orthogonal decomposition ($\mathbb{P}$ is the orthogonal projector on divergence-free vectorfields, and $\mathbb{Q}=I_d -\mathbb{P}$ is the orthogonal projector on gradients) and denoting $F_{\ro} \overset{def}{=}\div\left[ k(\ro) \n \ro \otimes \n \ro\right]$:
$$
\begin{cases}
\d_t u_L -\Delta u_L =-\kmu \mathbb{P} F_{\ro},\\
\n P_L =-\kmu \mathbb{Q} F_{\ro}
\end{cases}
$$
As the external force term is independant of $t$, we immediately obtain that:
\begin{equation}
\|\n P_L\|_{L_t^1 \dot{B}_{p,1}^{\frac{n}{p}-1}} \leq \kmu \|F_{\ro}\|_{L_t^1 \dot{B}_{p,1}^\frac{n}{p}} \leq \kmu C_{\ro} t.
\end{equation}
And concerning the velocity, following classical localization methods, for all $j\in \Z$ (see \cite{Dbook}):
$$
\ddj u_L (t) =e^{t\Delta} \frac{1}{\mub} \ddj u_0+ \kmu \int_0^t e^{(t-\tau)\Delta} \ddj F_{\ro} d\tau,
$$
taking the $L^p$-norm, thanks to the frequency localization and Lemma 2.4 from \cite{Dbook} and using that $F_{\ro}$ does not depend on $t$:
$$
\|\ddj u_L(t)\|_{L^p} \leq C\left(\frac{\|\ddj u_0\|_{L^p}}{\mub} e^{-ct2^{2j}} +\kmu \frac{1-e^{-ct 2^{2j}}}{c2^{2j}} \|\ddj F_{\ro}\|_{L^p} \right).
\label{refine}
$$
Then taking the $L_t^1$-norm and multiplying by $2^{j(\frac{n}{p}-1)}$ leads us to the classical refined estimate:
\begin{equation}
\|u_L\|_{L_t^1 \dot{B}_{p,1}^{\frac{n}{p}+1}} \leq C\left(\frac{1}{\mub} \Sum_{j\in \Z} (1-e^{-ct2^{2j}}) 2^{j(\frac{n}{p}-1)} \|\ddj u_0\|_{L^p} +\kmu t \|F_{\ro}\|_{\dot{B}_{p,1}^{\frac{n}{p}-1}} \right).
\end{equation}
Taking the $L_t^\infty$-norm in \eqref{refine} would lead to
$$
\|u_L\|_{L_t^\infty \dot{B}_{p,1}^{\frac{n}{p}-1}} \leq C\left(\frac{1}{\mub} \|u_0\|_{\dot{B}_{p,1}^{\frac{n}{p}-1}} +\kmu t \|F_{\ro}\|_{\dot{B}_{p,1}^{\frac{n}{p}-1}} \right),
$$
and taking the $L_t^2$-norm in \eqref{refine}:
\begin{equation}
\|u_L\|_{L_t^2 \dot{B}_{p,1}^\frac{n}{p}} \leq C\left(\frac{1}{\mub} \Sum_{j\in \Z} \sqrt{1-e^{-ct2^{2j}}} 2^{j(\frac{n}{p}-1)} \|\ddj u_0\|_{L^p} +\kmu \sqrt{t} \|F_{\ro}\|_{\dot{B}_{p,1}^{\frac{n}{p}-2}} \right).
\end{equation}
Thanks to the assumptions on $\ro$, we obtain that:
$$
\|F_{\ro}\|_{\dot{B}_{p,1}^{\frac{n}{p}-1} \cap \dot{B}_{p,1}^{\frac{n}{p}-2}} \leq (1+C(\|\ro\|_{L^\infty}) \|\ro-1\|_{\dot{B}_{p,1}^\frac{n}{p}}) \|\n \ro\|_{\dot{B}_{p,1}^\frac{n}{p}} \left(\|\n \ro\|_{\dot{B}_{p,1}^\frac{n}{p}}  +\|\ro-1\|_{\dot{B}_{p,1}^\frac{n}{p}} \right),
$$
and we end up with
\begin{equation}
\begin{cases}
\|(\d_t u_L, \n^2 u_L, \n P_L)\|_{L_t^1 \dot{B}_{p,1}^{\frac{n}{p}-1}} \leq C\left(\frac{1}{\mub} \Sum_{j\in \Z} (1-e^{-ct2^{2j}}) 2^{j(\frac{n}{p}-1)} \|\ddj u_0\|_{L^p} +\kmu t C_{\ro} \right),\\
\|u_L\|_{L_t^2 \dot{B}_{p,1}^\frac{n}{p}} \leq C\left(\frac{1}{\mub} \Sum_{j\in \Z} \sqrt{1-e^{-ct2^{2j}}} 2^{j(\frac{n}{p}-1)} \|\ddj u_0\|_{L^p} +\kmu \sqrt{t} C_{\ro} \right).
\end{cases}
\end{equation}
And to give an explicit example, under the stronger assumption $u_0\in \dot{B}_{p,1}^{\frac{n}{p}-1} \cap \dot{B}_{p,1}^{\frac{n}{p}+1}$, using that for all $\alpha \geq 0$, $1-e^{-\alpha} \leq \alpha$, the condition in \eqref{CondT2} is satisfied when:
$$
\begin{cases}
\Big(C \frac{\|u_0\|_{\dot{B}_{p,1}^{\frac{n}{p}+1}}}{\mub} +\kmu C_{\ro} \Big) t \leq \frac{\ee}{4},\\
\Big(C \frac{\|u_0\|_{\dot{B}_{p,1}^\frac{n}{p}}}{\mub} +\kmu C_{\ro}\Big) \sqrt{t} \leq \frac{\ee}{4},
\end{cases}
$$
Recalling that due to the time rescaling \eqref{rescale}, a lower bound for the lifespan corresponding to the original system is:
$$
T_0=\frac{\ee}{4} \min \left( \frac{1}{C \|u_0\|_{\dot{B}_{p,1}^{\frac{n}{p}+1}} +\kmu C_{\ro} \mub}, \frac{\ee}{4} \frac{\mub}{\Big(C \|u_0\|_{\dot{B}_{p,1}^\frac{n}{p}}+\kmu C_{\ro} \mub\Big)^2} \right).
$$
\subsection{Extension of the results}
It is usual in the Besov setting to try to extend the results in the case where the summation index $r$ is strictly greater than 1. First we recall that for the velocity, the case $r>1$ in the inhomogeneous Navier-Stokes system must be carefully studied as there may be difficulties to define the flow (as the velocity is not necessarily Lipschitz). If $r=1$ for the velocity, we may investigate the case $r>1$ for the density.

The question is then on one hand, to be able to estimate $\mu(\ro)-1$ and $k(\ro)-1$ in $\dot{B}_{p,r}^s$ (which requires $s<\frac{n}{p}$ or $s\leq \frac{n}{p}$ if $r=1$), and on the other hand, to estimate in $\dot{B}_{p,1}^{\frac{n}{p}-1}$ the capillary term:
\begin{multline}
\div\Big(k(\ro) \n \ro \otimes \n \ro\Big) =\Big(k'(1)+ (k'(\ro) -k'(1))\Big)\n \ro\cdot \n \ro \otimes \n \ro\\
+ \Big(k(1)+ (k(\ro)-k(1))\Big) \div(\n \ro \otimes \n \ro),
\end{multline}
where the last term is of the form $\n \ro \cdot \n^2 \ro$.

Estimating the capillary term as written in the left-hand side would require in particular to estimate in $\dot{B}_{p,1}^\frac{n}{p}$ the product $\n \ro \otimes \n \ro$ and due to the paraproduct $T_{\n \ro} \n \ro$ (see appendix) there are only two alternatives ($\|\n \ro\|_{L^\infty} \|\n \ro\|_{\dot{B}_{p,1}^\frac{n}{p}}$ or $\|\n \ro\|_{\dot{B}_{\infty,r}^{-\ee}} \|\n \ro\|_{\dot{B}_{p,\bar{r}}^{\frac{n}{p}+\ee}}$), both of them requiring that $\n \ro\in \dot{B}_{p,1}^\frac{n}{p}$ (as $\dot{B}_{p,\infty}^{\frac{n}{p}-\ee} \cap \dot{B}_{p,\infty}^{\frac{n}{p}+\ee} \hookrightarrow \dot{B}_{p,1}^\frac{n}{p}$).

In the right-hand side formulation, the same occurs for the last term due to the paraproduct $T_{\n \ro} \n^2 \ro$ so that we cannot relax the condition $r=1$ even in the constant coefficients case. In this case we easily adapt the result of \cite{RP2} to obtain local existence if $\n \ro \in \dot{B}_{p,1}^\frac{n}{p}$ and $\ro-1$ is small in the multiplier space $\mathcal{M}(\dot{B}_{p,1}^{\frac{n}{p}-1})$.

\section{Appendix}

The first part is devoted to a quick presentation of the Littlewood-Paley theory.

\subsection{Littlewood-Paley theory}

\label{appendiceLP}

In this section we briefly present the classical dyadic decomposition and some properties (for more details we refer to \cite{Dbook} chapter 2). Consider a smooth radial function $\chi$ supported in the ball $B(0, \frac{4}{3})$, equal to $1$ in a neighborhood of $B(0, \frac{3}{4})$ and such that $r\mapsto \chi(r.e_1)$ is nonincreasing over $\R_+$. If we define $\varphi(\xi)=\chi(\frac{\xi}{2})-\chi(\xi)$, $\varphi$ is supported in the annulus $\cC(0,\frac{3}{4}, \frac{8}{3})$ (equal to $1$ in a sub-annulus), and satisfy that for all $\xi\in \R^3\setminus \{0\}$,
$$
\sum_{q\in\Z} \varphi(2^{-q} \xi)=1.
$$
Then for all tempered distribution $u$ we define for all $q \in \mathbb{Z}$:
$$
\dot{\Delta}_q u=\mathcal{F}^{-1}\big(\varphi(2^{-q}\xi) \hat{u}(\xi)\big) \mbox{ and }\dot{S}_q u  =\sum_{p\leq q-1} \dot{\Delta}_p u= \chi(2^{-q}D) u.
$$
The homogeneous Besov spaces are defined as follows:
$$
\dot{B}_{p,r}^s=\{u \in \cS'(\R^3), \mbox{ with } \underset{q\rightarrow -\infty}{\mbox{lim}} \dot{S}_q u=0 \mbox{ and } \|u\|_{\dot{B}_{p,r}^s} \overset{\mbox{def}}{=} \|\big(2^{qs} \|\ddq u\|_{L^p} \big)_{q \in \Z}\|_{\ell^r} < \infty\}.
$$
\begin{rem}
\sl{Due to the support of $\varphi$, we easily obtain that
\begin{equation}
\dot{\D}_j \dot{\D}_l =0 \mbox{ if } |j-l|\geq 2.
\label{supports}
\end{equation}
}
\end{rem}
Let us now turn to the Bony decomposition, coming from the fact that for all distributions $u,v$, we can write (at least formally) the product as follows:
$$
uv=(\sum_{l\in \Z} \ddl u)(\sum_{j\in \Z} \ddj v).
$$
A more efficient way to write this product is the following Bony decomposition, where we basically set three parts according to the fact that the frequency $l$ of $u$ is of smaller, comparable or bigger size than the frequency $j$ of $v$:
\begin{equation}
uv= T_u v+ T_v u+ R(u, v),
\label{Bonydecp}
\end{equation}
where
\begin{itemize}
\item $T$ is the paraproduct~: $T_uv:=\sum_l\dot S_{l-1}u\ddl v$ (for each $l$, the term has its frequencies in an annulus of size $2^l$),
\item $R$ is the remainder~: $R(u, v)= \sum_l\sum_{|\alpha|\leq1}\ddl u\,\dot\Delta_{l+\alpha}v$ (the term has its frequencies in a ball of size $2^l$).
\end{itemize}
In this article we will often use the following estimates for the paraproduct and remainders in order to deal with nonlinear terms (we refer to \cite{Dbook} Section 2.6 for general statements, more properties of continuity for the paraproduct and remainder operators: 
\begin{prop}\label{p:op}
For any $(s,p,r)\in\R\times[1,\infty]^2$ and $t<0,$ there exists a constant $C$ such that 
$$
\|T_uv\|_{\dot B^s_{p,r}}\leq C\|u\|_{L^\infty}\|v\|_{\dot B^{s}_{p,r}}\quad\hbox{and}\quad
\|T_uv\|_{\dot B^{s+t}_{p,r}}\leq C\|u\|_{\dot B^t_{\infty,\infty}}\|v\|_{\dot B^{s}_{p,r}}.
$$
For any $(s_1,p_1,r_1)$ and $(s_2,p_2,r_2)$ in $\R\times[1,\infty]^2$ there exists a constant $C$
such that 
\begin{itemize}
\item if 
$s_1+s_2>0,$ $1/p:=1/p_1+1/p_2\leq1$ and $1/r:=1/r_1+1/r_2\leq1$ then
 $$\|R(u,v)\|_{\dot B^{s_1+s_2}_{p,r}}\leq C\|u\|_{\dot B^{s_1}_{p_1,r_1}}
\|v\|_{\dot B^{s_2}_{p_2,r_2}}.$$
\end{itemize}
\end{prop}
Let us now turn to the composition estimates. We refer for example to \cite{Dbook} (Theorem $2.59$, corollary $2.63$)):
\begin{prop}
\sl{\begin{enumerate}
 \item Let $s>0$, $u\in \dot{B}_{p,r}^s\cap L^{\infty}$ ($s<\frac{n}{p}$ and $r>1$ or $s\leq\frac{n}{p}$ of $r=1$) and $F$ a smooth function such that $F(0)=0$. Then $F(u)\in \dot{B}_{p,r}^s$ and there exists a function of one variable $C_0$ only depending on $s$, $n$ and $F'$ such that
$$
\|F(u)\|_{\dot{B}_{p,r}^s}\leq C_0(\|u\|_{L^\infty})\|u\|_{\dot{B}_{p,r}^s}.
$$
\item under the same assumptions, for any $u,v \in \dot{B}_{p,r}^s \cap L^\infty$, then there exists a constant $C$ such that
$$\|F(v)-F(u)\|_{\dot{B}_{p,r}^s} \leq C(k'', \|u\|_{\dot{B}_{p,r}^s \cap L^\infty}, \|v\|_{\dot{B}_{p,r}^s \cap L^\infty}) \left(\|v-u\|_{\dot{B}_{p,r}^s} +\|v-u\|_{L^\infty} \right)
$$
\end{enumerate}}
\label{estimcompo}
\end{prop}
\subsection{Lagrangian change of variables}
We gather in this section the main properties that we use in the process of the Lagrangian change of variable. For more details and proofs we refer to \cite{RP2} (the compressible version and more general results can be found in \cite{Dlagrangien2}).
\begin{prop}
\sl{
\label{Comp}Let $X$ be a globally defined bi-lipschitz diffeomorphism of $\mathbb{R}^{n}$ and $(s,p,q)$ with $1\leq p <\infty$ and $-\frac{n}{p^{\prime}}<s\leq\frac{n}{p}$. Then $a\rightarrow a\circ X$ is a self-map over $\dot{B}_{p,1}^{s}(\R^n)$ in the following cases:
\begin{enumerate}
\item $s\in (0,1)$,
\item $s\geq 1$ and $( DX-Id)  \in\dot{B}_{p,1}^\frac{n}{p}$.
\end{enumerate}}
\end{prop}

\begin{prop}
\sl{
\label{Formule}Let $X$ be a $C^{1}$ diffeomorphism over $\R^n$. For any function $f$ we introduce the notation $\bar{f}=f\circ X$. Then we obviously have ($Df$ denotes the jacobean matrix) :
\begin{equation}
\overline{Df}=D\bar{f} \cdot (DX)^{-1} \mbox{ and } \overline{\n f}=(\n X)^{-1}\cdot \n \bar{f}.
\end{equation}
Let $f$ be a $C^{1}$ scalar function over $\R^n$ and $F$ a $C^{1}$ vector field. If we assume in addition that the Jacobean determinant $J\overset{def}{=} \det\left(DX\right)$ is positive then we have  (we recall that for any matrix $A$, we introduce its adjugate matrix, denoted by $\mbox{adj}(A)$, i.-e. the transposed cofactor matrix. If $A$ is invertible $\mbox{adj}(A)=\mbox{det}(A) A^{-1}$):
\begin{align*}
\overline{\n f} & =J^{-1} \div(\mbox{adj} (DX) \bar{f}),\\
\overline{\div F} & =J^{-1} \div(\mbox{adj} (DX) \bar{F}),\\
\end{align*}
}
\end{prop}
As a consequence, for all vectorfield $F$:
\begin{multline}
J^{-1} \div(\mbox{adj} (DX) \bar{F})=\overline{\div F} =\overline{\mbox{tr} (DF)} =\mbox{tr}(\overline{DF})\\
=\mbox{tr}(D\bar{F}\cdot (DX)^{-1}) =D\bar{F} : (DX)^{-1},
\label{magic}
\end{multline}
which serves as a crucial ingredient in \cite{RP2, Dlagrangien2, CB}.

As stated in Remark \ref{ChgtLagr}, in order to perform the Lagrangian change of variable relatively to some $\bar{v}$ (a time dependent vector field) we first define its associated flow:
\begin{equation}
X_{\bar{v}}\left(  t,y\right)  =y+\int_{0}^{t}\bar{v}\left(\tau,y\right) d\tau,\label{flow}%
\end{equation}
and we denote $A_{\bar{v}}=\left(  DX_{\bar{v}}\right)^{-1}$, $\Jv =\mbox{det}(DX_{\bar{v}}^{-1})$.
\\

The aim is to rewrite the Navier-Stokes system (and \eqref{NSIK}) in a new form featuring a constant density and no transport terms. We refer to \cite{RP2} for the changes in the core terms of the inhomogeneous incompressible Navier-Stokes system. Let us just precise the computation for the capillary term. Thanks to Proposition \ref{Formule}:
\begin{multline}
\Big(\div [k(\rho) \n \rho \otimes \n \rho]\Big) \circ X_{\vb} =\Jv^{-1}\div \left[\adjv \Big((k(\rho) \n \rho \otimes \n \rho)\circ X_{\vb}\Big) \right]\\
=\Jv^{-1}\div \left[\adjv  k(\rhob) \Big(\n \rho \circ X_{\vb}\Big) \otimes \Big(\n \rho\circ X_{\vb}\Big) \right]\\
=\Jv^{-1}\div \left[k(\rhob) \adjv \Big({^t DX_{\vb}^{-1}} \n \rhob\Big) \otimes \Big({^t DX_{\vb}^{-1}} \n \rhob\Big) \right]\\
=\Jv^{-1}\div \left[k(\rhob) \adjv \Big({^t \Av} \n \rhob\Big) \otimes \Big({^t \Av} \n \rhob\Big) \right].
\end{multline}
And in the case where $X_{\vb}$ is measure-preserving ($\Jv\equiv 1$) we end up with:
$$
\Big(\div [k(\rho) \n \rho \otimes \n \rho]\Big) \circ X_{\vb} =\div \left[\Av^{-1} k(\rhob) \Big({^t \Av} \n \rhob\Big) \otimes \Big({^t \Av} \n \rhob\Big) \right].
$$
But we also need to go back to the original functions that is if $\ub$ is the solution of the transformed system we need to be able to define $u(t,.)=\ub(t,X_{\ub}^{-1}(t,.))$, this is why we need $X_{\vb}$ to be a global diffeomorphism, which is the object of the following proposition (we refer for example to the appendix of \cite{CB} for precisions and proofs):
\begin{prop}
\sl{
\label{PropozitieFinal}Let us consider $\bar{v}\in C_{b}(\left[  0,T\right]
,\dot{B}_{p,1}^{\frac{n}{p}-1})$ with $\partial_{t}\bar{v},$ $\nabla^{2}%
\bar{v}\in L_{T}^{1}(\dot{B}_{p,1}^{\frac{n}{p}-1})$. Then, there exists $\ee_0>0$ such that if
\begin{equation}
\int_0^T \left\Vert \nabla\bar{v}\right\Vert _{\dot{B}_{p,1}^{\frac{n}{p}}} dt \leq\ee_0, \label{smallnes12}%
\end{equation}
then, $X_{\bar{v}}$ introduced in $\left(  \text{\ref{flow}}\right)  $ is a global $C^{1}$-diffeormorphism over $\mathbb{R}^{n}$. Moreover, if we have $\div\left(  \adjv \bar{v}\right)=0$ then, $X_{\bar{v}}$ is measure preserving i.e. $\det DX_{\bar{v}}\equiv 1$.
}
\end{prop}
Under the smallness property \eqref{smallnes12}, the following properties (we refer to \cite{RP2, Dlagrangien2} for proofs) help us to deal with the various matrices introduced by the change of variables:
\begin{prop}
\sl{
Let us consider $\bar{v}\in E_{T}$ satisfying the smallness condition $\left(
\text{\ref{smallnes12}}\right)  $. Let $X_{v}$ be defined by $\left(
\text{\ref{flow}}\right)  $. Then for all $t\in\lbrack0,T]:$%
\begin{align}
\left\Vert Id-A_{\bar{v}}\left(  t\right)  \right\Vert _{\dot{B}_{p,1}^{\frac{n}{p}}} +\|I_d-\adjv(t)\|_{\dot{B}_{p,1}^{\frac{n}{p}}} &  \lesssim\left\Vert \nabla\bar{v}\right\Vert _{L_{t}^{1}\dot{B}_{p,1}^{\frac{n}{p}}},\label{A1}\\
\left\Vert \partial_{t}\adjv \left(  t\right)  \right\Vert _{\dot{B}_{p,1}^{\frac{n}{p}-1}}  &  \lesssim\left\Vert \nabla\bar{v}(t)\right\Vert_{\dot{B}_{p,1}^{\frac{n}{p}-1}},\label{A2}\\
\left\Vert \partial_{t}\adjv \left(  t\right)  \right\Vert _{\dot{B}_{p,1}^{\frac{n}{p}}}  &  \lesssim\left\Vert \nabla\bar{v}(t)\right\Vert_{\dot{B}_{p,1}^{\frac{n}{p}}},
\label{A3}\\
\|\Jv^{\pm}(t) -1\|_{\dot{B}_{p,1}^{\frac{n}{p}}} & \lesssim \|\n \vb\|_{L_t^1 \dot{B}_{p,1}^{\frac{n}{p}}}. \label{A3b}
\end{align}
}
\label{propestimL}
\end{prop}
\begin{prop}
\sl{
Let $\bar{v}_{1},\bar{v}_{2}\in E_{T}$ satisfying the smallness condition
$\left(  \text{\ref{smallnes12}}\right)  $ and $\delta v=\bar{v}_{2}-\bar
{v}_{1}$. Then we have:%
\begin{align}
\left\Vert A_{\bar{v}_{1}}-A_{\bar{v}_{2}}\right\Vert _{L_{T}^{\infty}\dot{B}_{p,1}^{\frac{n}{p}}} +\left\Vert \mbox{adj}(DX_{\bar{v}_{1}})-\mbox{adj}(DX_{\bar{v}_{2}})\right\Vert _{L_{T}^{\infty}\dot{B}_{p,1}^{\frac{n}{p}}} &  \lesssim\left\Vert \nabla\delta v\right\Vert_{L_{T}^{1}\dot{B}_{p,1}^{\frac{n}{p}}},\label{A4}\\
\left\Vert \d_t \mbox{adj}(DX_{\bar{v}_{1}})-\d_t \mbox{adj}(DX_{\bar{v}_{2}})\right\Vert_{L_{T}^{1}\dot{B}_{p,1}^{\frac{n}{p}}} +\|J_{\bar{v}_1}^\pm(t) -J_{\bar{v}_2}^\pm(t)\|_{\dot{B}_{p,1}^{\frac{n}{p}}} &  \lesssim\left\Vert \nabla\delta v\right\Vert _{L_{T}^{1}\dot{B}_{p,1}^{\frac{n}{p}}},\label{A5}\\
\left\Vert \d_t \mbox{adj}(DX_{\bar{v}_{1}})-\d_t \mbox{adj}(DX_{\bar{v}_{2}})\right\Vert_{L_{T}^{2} \dot{B}_{p,1}^{\frac{n}{p}-1}}  &  \lesssim\left\Vert \nabla\delta v\right\Vert _{L_{T}^{2}\dot{B}_{p,1}^{\frac{n}{p}-1}}, \label{A6}
\end{align}
}
\label{propestimL2}
\end{prop}

\textbf{Aknowledgements :} This work was supported by the ANR project INFAMIE, ANR-15-CE40-0011. The first author is partially supported by a grant of the Romanian National Authority for Scientific Research and Innovation, CNCS–UEFISCDI, project number PN-II-RU-TE-2014-4-0320.


\begin{thebibliography}{}

\bibitem{Abels} Abels, H. Strong well-posedness of a diffuse interface model for a viscous, quasi-incompressible two-phase flow, \textit{Siam J. Math. Anal}, Vol 44 (1) (2012), p 316-340.





\bibitem{Aki} Aki, G., Dreyer, W., Giesselmann, J., Kraus, C. A quasi-incompressible diffuse interface model with phase transition, \textit{Math. Models Methods Appl. Sci.} 24 (2014), no. 5, 827-861.

\bibitem{Dbook} Bahouri, H., Chemin, J.-Y., Danchin, R. Fourier analysis and nonlinear partial differential equations, \textit{Grundlehren der mathematischen Wissenschaften}, \textit{343}, \textit{Springer Verlag}, 2011.


\bibitem{BDL} Bresch, D., Desjardins, B. and Lin, C.-K. On some compressible fluid models: Korteweg,lubrication and shallow water systems. \textit{Comm. Partial Differential Equations}, {\bf 28(3-4)} : 843-868, 2003.

\bibitem{B2} Bresch, D., Guill\'en-Gonz\'alez, F., Lemoine, J., Rodr\'iguez-Bellido,  M.-A., Sy, M. Local strong solution for the incompressible Korteweg model, \textit{C. R. Acad. Sci. Paris}, Ser. I 342 (2006), 169-174.

\bibitem{CB} Burtea, C. Optimal well-posedness for the inhomogeneous incompressible Navier-Stokes system with general viscosity, \textit{submitted}, 
https://hal.archives-ouvertes.fr/hal-01349093


\bibitem{CH1} Charve, F., Haspot, B. Convergence of capillary fluid models: from the non-local to the local Korteweg model, \textit{Indiana U. Math. J.}, \textbf{60}(6) (2011), 2021-2060.


\bibitem{CH4} Charve, F., Haspot, B. On a Lagrangian method for the convergence from a non-local to a local Korteweg capillary fluid model, \textit{J. Funct. Anal.} \textbf{265} (7) (2013), p.1264-1323.

\bibitem{FCorder} Charve, F. Convergence of a low order non-local Navier-Stokes-Korteweg system: the order-parameter model, to appear in \textit{Asymptotic Analysis}.

\bibitem{FCLarge} Charve, F. Local in time results for local and non-local capillary Navier-Stokes systems with large data, \textit{J.  of Diff.  Eq.} \textbf{256} (7) (2014), pp 2152-2193.


\bibitem{Noboru} Chikami, N., Ogawa, T. On the well-posedness of the incompressible Navier-Stokes-Poisson system in Besov spaces, to appear in J. Evol. Eq.

\bibitem{DD} Danchin, R., Desjardins, B. Existence of solutions for compressible fluid models of Korteweg type, \textit{Annales de l'IHP, Analyse Non Lin\'eaire}, {\bf 18} (2001), 97-133.



\bibitem{Dtruly} Danchin, R. Well-Posedness in Critical Spaces for Barotropic Viscous fluids with truly not constant density, \textit{Comm. in Partial Diff. Eq.}, \textbf{32} (2007), 1373--1397.

\bibitem{Dlagrangien} Danchin, R. Uniform estimates for transport-diffusion equations, \textit{J. Hyp. Diff. Eq.}, \textbf{4}(1), 1-17 (2007).



\bibitem{RP2} Danchin, R. and Mucha, P. B. (2012). A Lagrangian approach for the incompressible Navier-Stokes equations with variable density. \textit{Comm. Pure Appl. Math.}, 65(10):1458-1480.


\bibitem{Dlagrangien2} Danchin, R. (2014). A Lagrangian approach for the compressible Navier-Stokes equations. \textit{Ann. Inst. Fourier (Grenoble)}, 64(2):753-791.

\bibitem{RPbook} Danchin, R. and Mucha, P. B. (2015). Critical functional framework and maximal regularity in action on systems of incompressible flows. \textit{Mém. Soc. Math. Fr. (N.S.)}, (143):vi+151.

\bibitem{Has1} Haspot, B. Existence of strong solutions for nonisothermal Korteweg system, \textit{Annales Math\'ematiques Blaise Pascal} {\bf 16}, 431-481 (2009).

\bibitem{Has2} Haspot, B. Cauchy problem for capillarity Van der Waals model, \textit{Hyperbolic problems: theory, numerics and applications}, 625634, Proc. Sympos. Appl. Math., 67, Part 2, Amer. Math. Soc., Providence, RI, 2009.


\bibitem{Has4} Haspot, B. Existence of global strong solution for the compressible Navier-Stokes system and the Korteweg system in two-dimension,  \textit{Methods and Applications of Analysis}, Vol. 20 (2013), No. 2, pp. 141–164.

\bibitem{Has5} Haspot, B. Existence of global strong solution for Korteweg system with large infinite energy initial data, \textit{Journal of Mathematical Analysis and Applications}, 438 (2016), no. 1, 395–443



\bibitem{TH1} Hmidi, T. R\'egularit\'e h\"olderienne des poches de tourbillon visqueuses, \textit{Journal de Math\'ematiques pures et appliqu\'ees}, \textbf{84}(11), (2005) 1455-1495.




\bibitem{PZ1} Paicu, M. and Zhang, P. (2012). Global solution to the 3d incompressible inhomogeneous Navier-Stokes system. \textit{J. Funct. Anal.}, 262:3556-3584.

\bibitem{PZ2} Paicu, M., Zhang, P., and Zhang, Z. (2013). Global unique solvability of inhomogeneous Navier-Stokes equations with bounded density. \textit{Communications in Partial Differential Equations}, 38(7):1208-1234.


\bibitem{Ro2} Rohde C., On local and non-local Navier-Stokes-Korteweg systems for liquid-vapour phase transitions. \textit{ZAMM Z. Angew. Math. Mech.} {\bf 85} (2005), no. 12, 839-857.

\bibitem{Rohdeorder} Rohde C., A local and low-order Navier-Stokes-Korteweg system, \textit{Nonlinear partial differential equations and hyperbolic wave phenomena, 315-337, Contemp. Math.}, \textbf{526 (2010)}, \textit{Amer. Math. Soc., Providence, RI.}


\bibitem{YYL} Yang, J., Yao, L., Zhu, C., Vanishing capillarity–viscosity limit for the incompressible inhomogeneous fluid models of Korteweg type, \textit{Z. Angew. Math. Phys.} (2015) 66(5), 2285-2303. 

\end{thebibliography}
\end{document}